\documentclass[12pt]{article}

\usepackage{amsfonts}
\usepackage{amsfonts}  
\usepackage{amssymb}
\usepackage{amsthm}
\usepackage{amsmath}
\usepackage{color}

\title{Qualitative properties of   a nonlinear system involving the $p$-Laplacian operator.} 

\author{F. Demengel}

\date{}

\date{}

\catcode`@=11 \@addtoreset{equation}{section} \catcode`@=12

\newtheorem{theo}{Theorem}[section]
\newtheorem{prop}[theo]{Proposition}
\newtheorem{rema}[theo]{Remark}

\newtheorem{lemme}[theo]{Lemma}

\def\R{\mathbb  R}

\begin{document}

\maketitle
\maketitle
\begin{abstract}
 In this article we consider  the nonlinear system involving the $p$-Laplacian 
 $$\left\{\begin{array}{lc}
 |u^\prime |^{p-2} u^{\prime \prime}  = u^{p-1} v^p&\\
 |v^\prime |^{p-2} v^{\prime \prime}  = v^{p-1} u^p&\ {\rm on} \ \R, \\
  u\geq 0,  v\geq 0&
\end{array}\right.$$
 for which we prove  symmetry, asymptotic behavior and  non degeneracy properties.  This  can help to a better understanding to what happens in the $N$ dimensional case, for which several authors prove a De Giorgi Type result under some additional growth and monotonicity assumptions. 
 \end{abstract}
 
\section{Introduction}

In this article we extend some of the results obtained in \cite{BLWZ} in the case  of the Laplacian, to the $p$-Laplacian case. More precisely we consider the system in $\R^N$ :

$$\left\{\begin{array}{lc}
{\rm div} (|\nabla u|^{p-2} \nabla u) &= u^{p-1} v^p\\
{\rm div} (|\nabla v|^{p-2} \nabla v) &= v^{p-1} u^p, 
\end{array}\right.$$
where $u$ and $v$ are supposed to be positive. 

In the  case $p=2$ this problem   comes from a phase  separation model.  As an example the Gross Pitaevskii, \cite{CLL} model  can be described by the  non linear elliptic system  

\begin{equation}\label{eqmodel}\left\{ \begin{array}{lc}-\Delta u + \alpha u^3 + \Lambda v^2 u = \lambda_{1, \Lambda} u \ & {\rm in} \ \Omega \\
-\Delta v + \alpha v^3 + \Lambda u^2 v = \lambda_{2, \Lambda} v \ & {\rm in} \ \Omega\\
 u>0, \ v>0, {\rm in} \ \Omega,    u = v = 0 \ &{\rm on}\  \partial \Omega \ \\
  \int_\Omega u^2 = \int_\Omega v^2 = 1& 
  \end{array}\right.
  \end{equation}
  where $\Omega$ is a smooth bounded domain in $\R^N$ and $\alpha, \beta$ are positive parameters, $\Lambda$ will become large.  
   Assuming that  $\sup ( \lambda_{1, \Lambda} ,  \lambda_{2, \Lambda} )\leq C$ for some constant independent of $\Lambda$, formally and up to subsequences $(u_\Lambda, v_\Lambda)$ converges to some pair $(u,v)$ which satisfies  $uv = 0$ and the equations
   \begin{equation}
 \label{eqlimit}\left\{ \begin{array}{lc}
-\Delta u + \alpha u^3  = \lambda_{1, \Lambda} u \ & {\rm in} \ \Omega_u = \{ x, u(x) >0\}\\
-\Delta v + \beta v^3  = \lambda_{2, \Lambda} v \ & {\rm in} \ \Omega_v = \{ x, v(x) >0\}\\
  \end{array}\right.
  \end{equation}
Several papers treat the convergence of
 $(u_\Lambda, v_\Lambda)$ away the interface $\gamma = \{ x, u(x) =0 = v(x)\}$, see for example  \cite{WW}  and \cite{CTV} , \cite{NTTV} for the uniform equicontinuity of $(u_\Lambda, v_\Lambda)$. 
 
 Near the interface, the profile of bounded solutions  of (\ref{eqmodel}) of the blow up equation is a system, 
 which is completely classified in the one dimensional case, \cite{BLWZ}, \cite{BTWW}.  This system is the following 
 $$\left\{ \begin{array}{lc}
  U^{\prime\prime} = U V^2& \ \\
   \ V^{\prime\prime} = V U^2 & \ {\rm in} \ \R, \\
 U>0,  V>0.&
 \end{array}
 \right.$$
 
In \cite{BLWZ}, the authors   expect  that the same system occurs in the $N$ dimensional case, say 
 
 $$\left\{ \begin{array}{lc}
  \Delta U = U V^2& \  \\
   \ \Delta V = V U^2 & \ {\rm in} \ \R^N, \\
 U>0,  V>0. &
 \end{array}
 \right.$$

 Furthermore they conjecture  that for any dimension $N\leq 8$,  the  system is  in fact one dimensional. They obtain this result  under some additional assumption on the growth of the solution, in dimension $2$.  The assumption $N\leq 8$ is  motivated by the case of scalar equations for which it is known that the scalar equation is not  necessary one dimensional
  for $N \geq 9$, \cite{DKW}.  The condition on the growth in the two dimensional case is satisfied in particular if  the solution of the system
   is at most linear at infinity,   as it is proved   in the case  $N=1$. When $N >2$, this is not sufficient, and  up to now, even 
   in the case $N=2$, this increasing behaviour  is not proved. 
   
    In \cite{F2} the author improved the result by establishing that  for $N=2$, as soon as $u$ and $v$ have at most algebraic increasing    behavior and satisfy for some component $\partial_N u >0$, $\partial_N v <0$,  then the solution is one dimensional. 
     This result is recently improved by Farina and Soave  by replacing the condition $\partial_N u >0$, $\partial_N v <0$ by the weaker condition 
     $\lim_{x_N\rightarrow \pm \infty} u(x^\prime, x_N)- v(x^\prime, x_N) = \pm \infty$, uniformly with respect to $x^\prime$, conserving the  assumption of algebraic growth. 
     Let us cite also the recent result of K. Wang \cite{W} which replaces the  monotonicity condition by the fact that $(u,v)$ is a local minimizer.

The present paper is motivated by the asymptotic study of  ( \ref{eqmodel}),  say 
 \begin{equation}
 \label{eqmodelp}\left\{ \begin{array}{lc}
-{\rm div} (|\nabla u|^{p-2} \nabla  u )+ \alpha u^{p+1} + \Lambda v^p u^{p-1} = \lambda_{1, \Lambda} u^{p-1} \  {\rm in} \ \Omega &\ \\
-{\rm div} (|\nabla v|^{p-2} \nabla  v )+ \beta v^{p+1} + \Lambda u^p v^{p-1} = \lambda_{2, \Lambda} v^{p-1} \  {\rm in} \ \Omega\ & \\
 u>0, \ v>0 {\rm in} \ \Omega,    u = v = 0 \ {\rm on}\  \partial \Omega \ & \\
  \int_\Omega u^p = \int_\Omega v^p = 1,  &  
  \end{array}\right.
  \end{equation}
  where $\Omega$ is a smooth bounded domain in $\R^N$ and $\alpha, \beta$ are positive parameters, $\Lambda$ will become large.  
   Such a pair of solutions is a critical point for the functional
  
  $E_\Lambda (u, v) = {1\over p} \int_\Omega( |\nabla u|^p + |\nabla v|^p)+ {\alpha \over p+2} \int_\Omega  |u|^{p+2} + {\beta \over p+2} \int_\Omega |v|^{p+2} + {\Lambda \over p} \int_\Omega |u|^p |v|^p$ under the constraint 
  $\int_\Omega |u|^p = \int_\Omega |v|^p = 1$.  
  
  Assume that  there exists some constant $C$   independent on $\Lambda$ with 
   
   $\sup_\Lambda (\lambda_{1, \Lambda}, \lambda_{2, \Lambda})\leq C$, for $\Lambda$ large . As $\Lambda$ goes to infinity,
    and up to subsequence $(u_\Lambda, v_\Lambda)$  tends formally to some pair $(u, v) $ which satisfies

      \begin{equation}\label{eq3}\left\{ \begin{array}{lc}
 -{\rm div} (\nabla u|^{p-2}\nabla u) + \alpha u^{p+1}  = \lambda_{1, \Lambda}  u^{p-1}& \ {\rm in} \ \Omega_u \\
  -{\rm div} (|\nabla v|^{p-2} \nabla v ) + \beta v^{p+1} = \lambda_{2, \Lambda} v^{p-1}& {\rm in} \ \Omega_v 
  \end{array}
   \right.
   \end{equation}
   where 
   $\Omega_u = \{ x, u(x) >0\}$,  $\Omega_v = \{ x, v(x) >0\}$. 
   
   It is not  our purpose here  to  follow this way. We  are interested in the one dimensional case and especialy 
     in the behavior of the  limit pair of solutions $(u, v)$  near the interface $\gamma = \{ x\in \Omega, u(x) = v(x) = 0\}$. 
    
  When $N=1$ and $\Omega = ]a, b[$ one has the result
     
     \begin{theo}\label{thlambda}
     Assume that $\Omega = ]a, b[$,  that 
     $(u_\Lambda , v_\Lambda)
$ solves  (\ref{eqmodelp}). There exists  $x_\Lambda \in \Omega $ such that 
     $u_\Lambda (x_\Lambda ) = m_\Lambda= v(x_\lambda) $ goes to zero, and $x_\Lambda$ tends to some point   in $ \gamma $. Furthermore if 
$\tilde u_\Lambda ={1\over m_\Lambda}  u_\Lambda (m_\Lambda y+ x_\Lambda)$ , $\tilde v_\Lambda ={1\over m_\Lambda}  v_\Lambda (m_\Lambda y+ x_\Lambda)$  and $y\in ]{a-x_\Lambda \over m_\lambda}, {b-x_\Lambda \over m_\Lambda}[$, $(\tilde u_\Lambda, \tilde v_\Lambda) $ converges locally uniformly to some pair $U, V$ which satisfies 
\begin{equation}\label{eq4}
\left\{\begin{array}{lc}
( |U^\prime |^{p-2} U^\prime )^\prime = V^p U^{p-1} \  &, \\
 ( |V^\prime |^{p-2} V^\prime )^\prime = U^p V^{p-1} \ & {\rm on} \ \R\\
U(0) = V(0) = 1. &  U, V >0.
\end{array}\right.
\end{equation}

Furthermore there exists  some  positive constant $T_\infty $ such that 
$${|U^\prime|^p } + {|V^\prime |^p}  -{U^p V^p } = {T_\infty}.  $$
\end{theo}

Next we are interested in the existence and  in the properties of the solutions of (\ref{eq4}). The existence of a non trivial solution is given in Theorem \ref{theoexi4}. In a second time  we  prove the  following 

 \begin{theo}\label{thUV}
  Let $N = 1$ and $(U$ $V)$ be  a non negative solution of 
  \begin{equation}\label{eq5}
  \left\{\begin{array}{lc}
  ( |U^\prime |^{p-2} U^\prime )^\prime = V^p U^{p-1} \ & , \\
 ( |V^\prime |^{p-2} V^\prime )^\prime = U^p V^{p-1} \ & {\rm on} \ \R.\\
\end{array}\right.
\end{equation}
 Then  up to exchanging $U$ and $V$, 
 
 1)   Up to translation $V(y) = U(-y)$. 
   
    2) $U^\prime >0$ everywhere, $U^\prime (+\infty) = (T_\infty )^{1\over p} $,  and  
    there exist some positive constants, $m, M$, $k$, $K$, $c$, $C$   such that near $-\infty$,  $m e^{-K x^2}\leq U(x) \leq M e^{-k x^2}$,  $c|x| U\leq U^\prime \leq C|x| U$.  
    
     Symmetric estimates  hold for $V$ exchanging $-\infty$ and $+\infty$.

3) Suppose that $\phi, \psi$ is a bounded solution of the linearized system
$$\left\{\begin{array}{lc}
 (|U^\prime |^{p-2} \phi^\prime)^\prime = (p-1)U^{p-2} V^p \phi + p V^{p-1} U^{p-1} \psi, \  & \\
     (|V^\prime |^{p-2} \psi^\prime)^\prime = (p-1)V^{p-2} U^p \psi + p U^{p-1} V^{p-1} \phi& 
     \end{array}\right.
     $$
in $\R$, then  there exists  some constant $c$ such that 
$(\phi, \psi) = c (U^\prime, V^\prime)$.

\end{theo}

We end this introduction by some reflexions about the De Giorgi's conjecture,  which , even if we do not treat it here,  is  after all,  at the origin of the present paper.

 As we said in the abstract,   the previous classification is a first step  if one want to prove a De Giorgi type result on the system 
 $$\left\{ \begin{array}{lc}
   {\rm div} (|\nabla u|^{p-2} \nabla u) = u^{p-1} v^p,&  \\
    {\rm div} (|\nabla v|^{p-2} \nabla v) = u^{p} v^{p-1} &{\rm in} \ \R^N\ 
    \end{array}\right.$$
    
    In \cite{DP}, \cite{DPP} the authors consider a more general system than the present one,  the quasilinear operator she studies includes the $p$Laplacian operator, and the right hand side included the case  studied here. She proves that under some  condition of growth on the solution, together with some  stability assumption on the couple of solutions,  then the solution is in fact one dimensional. The stability is in particular implied when the solution is "half monotone ", i.e.  if there exists one direction say $e_1$ such that 
    $\partial_1 u>0$ and $\partial_1 v <0$.  It is an exercise to prove that the growth condition (1.14) in \cite{DP} is  
    satisfied when the solution has at most linear growth at infinity,  in dimension 2, hence we recover a generalization of the  results in \cite{BTWW} in  the case $p\neq 2$.    

    Several questions are of interest: As we saied above,  in   \cite{F2}, and \cite{FS} the authors suppose that the solution has  at most  algebraic growth.  Can we have the same result in the $p$Laplacian case? 
      The answer is not immediate, since  the proof of Farina  and Soave   relies on some properties of  the Almgren frequency function which do not extend to the $p$Laplacian case .

     Another  interesting question is the following : Suppose that one replaces the Laplacian  system by a Fully  Nonlinear system. 
       Of course for Pucci's operators the one dimensional system is reduced, up to constant,  to the  Laplacian case.   But due to the non differentiability of the Pucci's operators, the definition of  stable solutions must be precised.         
         In the same order of ideas, one can imagine to treat the case of Fully Nonlinear degenerate or singular systems,   based on the model of the $p$Laplacian type treated here,  but not under divergence form, as  the following 
         
         $$\left\{ \begin{array}{lc}
         |\nabla u|^{\alpha} F(D^2 u) = u^{\alpha +1} v^{\alpha+2}& \\
           |\nabla v|^{\alpha} F(D^2 v) = v^{\alpha +1} u^{\alpha+2}& {\rm in} \ \R^N
           \end{array}\right.$$
            where $\alpha $ is some number $>-1$ and $F$ is fully nonlinear elliptic. The reader may consult  \cite{BD1} for properties of such operators and  a convenient definition of viscosity solutions.

          It would be far too long to cite  all the papers written about the  De Giorgi type result in the case of one  equation in place of a system. 
          Let us cite for the $p$-Laplacian case  the recent paper of  A. Farina and Valcinocci \cite{FV1}, and the very complete paper of Savin et al. \cite{SSV}.   
          For variations on the subject on De Giorgi's conjecture in  the case of a single equation the reader may consult  \cite{DeG}, \cite{F1}, \cite{FV2}, \cite{FSV}, \cite{BBG}, \cite{BHM}, \cite{GG}, \cite{AC}, \cite{Sa}, \cite{BD2}.        
    
     The paper is organized as follows :   In Section 2, we consider the one dimensional system 
     defining $(u_\Lambda, v_\Lambda) $ and prove  Theorem \ref{thlambda}.  Section 3 is devoted to the proof of Theorem \ref{thUV}. 
      Most of the technical details  of this section are postponed to the appendix in section 4.
      
     \bigskip
    
     {\it Acknowledgment}
     The author wishes to thank Alberno Farina and Enrico Valdinoci for having pointed  out several papers on the question.

\section{The original problem : Proof of theorem \ref{thlambda}}

In all that section we will frequently use in place of sequences, subsequences, without mentioning  it. 
 Let us consider for $\alpha$, $\beta$  and $\Lambda$, $\lambda_{1, \Lambda},\lambda_{2, \Lambda}$  some given positive constants 
 \begin{equation}\label{eq0}\left\{ \begin{array}{lc}
 -|u_\Lambda^\prime |^{p-2} u_\Lambda^{\prime\prime} + \alpha u_\Lambda^{p+1} + \Lambda v_\Lambda^p u_\Lambda^{p-1} = \lambda_{1, \Lambda} u_\Lambda^{p-1}& \ \\
  -|v_\Lambda^\prime |^{p-2} v_\Lambda^{\prime\prime} + \beta u_\Lambda^{p+1} + \Lambda v_\Lambda^{p-1} u_\Lambda^{p} = \lambda_{2, \Lambda} v_\Lambda^{p-1}&{\rm in} \ ]a, b[. \\
   u, \ v >0, u(a) = u(b)= v(a)= v(b)=0, & \int _a^b |u|^p = \int_a^b |v|^p = 1, 
   \end{array}
   \right.
   \end{equation}
   
    $(u_\Lambda,v_\Lambda)$ is then a solution of the minimizing eigenvalue problem
    
    $$\inf_{|u|_p = |v|_p= 1, (u, v) \in W^{1,p} (]a, b[)}   E_\Lambda (u,v) , $$
     where 
     
     $$E_\Lambda (u,v) = {1\over p} \int_a^b|u^\prime |^p + {1\over p} \int_a^b |v^\prime |^p + \alpha \int_a^b {|u|^{p+2}\over p+2} +  \beta \int_a^b {|v|^{p+2}\over p+2} +{\Lambda\over p} \int_a^b|u^p |v|^p.  $$
     
      Assume that 
      $$ \max_\Lambda (\lambda_{1, \Lambda}, \lambda_{2, \lambda}) \leq C, $$

        Due to the first equation in   (\ref{eq0})  multiplied by $u_\Lambda $,  and integrated  over $]a, b[$, one gets 
         \begin{equation}\label{eqlambda1}
         \int_a^b |u^\prime_\Lambda|^p + \alpha \int_a^b u_\Lambda^{p+2} +\Lambda \int_a ^b u_\Lambda^p v_\Lambda^p  = \lambda_{1, \Lambda}.
         \end{equation}
          Doing the same for $v_\Lambda$, one gets that $(u_\Lambda, v_\Lambda) $ is bounded in $W^{1,p}(]a, b[)^2$. 
          
          Furthermore
          
           \begin{lemme}\label{lemLambda}
           
            There exist $T_\Lambda$, $C_1$ and $C_2$ independent on $\Lambda$ such that 
            
             \begin{equation}\label{Tlambda}{|u^\prime_\Lambda |^p\over p} + { |v^\prime_\Lambda  |^p \over p} -\Lambda {u_\Lambda ^p v_\Lambda ^p \over p} -\alpha {u_\Lambda ^{p+2} \over p+2} -\beta {v_\Lambda ^{p+2}\over p+2} + \lambda_{1, \Lambda} {u_\Lambda ^p\over p} + \lambda_{2, \Lambda} {v_\Lambda ^p\over p} = {T_\Lambda\over p} 
             \end{equation}
             and 
             $$0< C_1 \leq T_\Lambda \leq C_2 < \infty$$
             \end{lemme}

             Proof 
             
              Multiply the first equation in (\ref{eq0}) by $u_\Lambda^\prime$,  the second one by  $v_\Lambda^\prime $,     and add the two equations, we obtain that 
              
              $$ -({|u^\prime_\Lambda|^p\over p} )^\prime -({|v^\prime_\Lambda|^p\over p} )^\prime+ 
              \alpha {(u_\Lambda^{p+2})^\prime \over p+2}  + \beta {(v_\Lambda^{p+2})^\prime \over p+2}  + \Lambda{ (v_\Lambda^p u_\Lambda^p)^\prime \over p} = \lambda_{1, \Lambda} {(u_\Lambda^p)^\prime \over p}+ \lambda_{2, \Lambda}  {(v_\Lambda^p)^\prime \over p}$$
               Hence there exists some constant $T_\Lambda$ such that (\ref{Tlambda}) is satisfied.  Integrating on $[a,b]$  the equation defining $T_\Lambda$, one  gets 
               
                $$\int_a^b{|u_\Lambda^\prime |^p\over p} + \int_a^b{ |v_\Lambda^\prime |^p \over p} -\Lambda \int_a^b{u_\Lambda^p v_\Lambda^p \over p} -\alpha \int_a^b{u_\Lambda^{p+2} \over p+2} -\beta\int_a^b {v_\Lambda^{p+2}\over p+2} + \lambda_{1, \Lambda}  + \lambda_{2, \Lambda}  ={ T_\Lambda(b-a)\over p}.$$

                On the other hand, using (\ref{eqlambda1}) for $u_\Lambda$ and  its analogous for $v_\lambda$, one has   
                 $$\int_a^b |u^\prime_\Lambda|^p + \int_a^b |v^\prime_\Lambda|^p+ \alpha \int_a^b u_\Lambda^{p+2} +\beta \int_a^b v_\Lambda^{p+2} 
                 +2\Lambda \int_a^b u_\Lambda^p v_\Lambda^p      = \lambda_{1, \Lambda} +\lambda_{2, \Lambda} . $$
            Combining the two equations one gets 
                  
                  $$T_\Lambda (b-a) = 2 \int_a^b |u^\prime_\Lambda|^p + 2 \int_a^b |v^\prime_\Lambda|^p + { 2\alpha\over p+2}\int_a^b u_\Lambda^{p+2} + {2\beta \over p+2}  \int_a^b v_\Lambda^{p+2}+\Lambda \int_a^b u_\Lambda^p v_\Lambda^p  $$
                  and using Poincar\'e's inequality
                  $\int_a^b |u^\prime _\Lambda|^p \geq C \int_a^b |u_\Lambda|^p \geq C$, 
                 one obtains  that 
                 $T_\Lambda\geq C>0$.

                   On the other hand since $(u_\Lambda)$ and $(v_\Lambda)$ 
                   are bounded in $W^{1,p}(]a, b[)$ and since $\lambda_{1,\Lambda}$ and $\lambda_{2, \Lambda}$ are bounded, one gets that 
                   $T_\Lambda$ is bounded from above.

                   Furthermore using  the equation defining $T_\Lambda$,  and the fact that $u_\Lambda$ and $v_\Lambda$ 
                   vanish  on the end points, one gets that 
                $u^\prime_\Lambda (a)$ , $u^\prime_\Lambda (b)$, $v^\prime_\Lambda (a)$, $v^\prime_\Lambda (b)$ 
                    are bounded independently on $\Lambda$. Integrating  the first equation in (\ref{eq0}) between $a$ and $x$ one gets 
                   $$|u_\Lambda^\prime |^{p-2} u_\Lambda^\prime (x) -|u_\Lambda^\prime |^{p-2} u_\Lambda^\prime (a) = \alpha 
                   \int_a^x u_\Lambda^{p+1} + \Lambda \int_
                   a^x v_\Lambda ^p u_\Lambda^{p-1} + \lambda_{1, \Lambda}  \int_a^x u_\Lambda^{p-1}$$
                    using this on $\{x=b\}$ one gets  $\Lambda \int_a^b u_\Lambda^{p-1} v_\Lambda^p \leq C$,  and by the positivity  that 
                                        $|u_\Lambda^\prime |^{p-2} u_\Lambda^\prime (x)-|u_\Lambda^\prime |^{p-2} u_\Lambda^\prime (a) \leq C$.
                     Doing the same between $x$ and $b$ one obtains that 
                     $|u_\Lambda^\prime |\leq C$
                      for  some constant independent on $\Lambda$. In the  same  manner,  $|v_\Lambda^\prime|\leq C$. 
                  
                    \begin{theo}\label{th1} 
                    
                    Assume that $ (u_\Lambda, v_\Lambda)$ solves the system (\ref{eq0}) in $]a, b[$. There exists $x_\Lambda$ such that 
                     $m_\Lambda = u(x_\Lambda ) = v_\Lambda (x_\Lambda) \rightarrow 0$ as $\Lambda$ goes to infinity,   and 
                    $x_\Lambda \rightarrow x_\infty \in ]a, b[$. Furthermore
                    $m_\Lambda^{2p} \Lambda \rightarrow C_o>0$ and  $\Lambda^{1\over 2p} \min (x_\Lambda -a, b-x_\Lambda) \rightarrow +\infty$.
                    \end{theo}

Proof of Theorem \ref{th1}

 By the previous estimates   $(u_\Lambda)$ and $(v_\Lambda)$ are  relatively compact in ${\cal C}([a, b])$. In particular up to subsequence,  
 $u_\Lambda$ and $v_\Lambda$ are uniformly convergent.  Let $(u_\infty, v_\infty)$ be the limit of  such subsequence. By the identity $\int_a^b |u_\infty|^p = \int_a^b |v_\infty|^p$, and by the uniform convergence ,   there exists $x_\infty$ and  $x_\Lambda$ which tends to $x_\infty$,  such that $m_\Lambda = u_\Lambda(x_\Lambda) = v_\Lambda (x_\Lambda)$ and 
 $x_\infty\in \{ u_\infty(x) = v_\infty(x)\} = \gamma$.

                    To prove that
                      
                      $$\limsup \Lambda m_\Lambda ^{2p} < \infty,$$
                       we  argue by contradiction and define 
                      $\tilde u_\Lambda = {1\over m_\Lambda} u({y\over m_\Lambda \Lambda^{1\over p}} + x_\Lambda)$,  $\tilde v_\Lambda = {1\over m_\Lambda} v({y\over m_\Lambda \Lambda^{1\over p}} + x_\Lambda)$
                        where $y\in ]-x_\Lambda + a m_\Lambda \Lambda^{1\over p}, -x_\Lambda + b m_\Lambda \Lambda^{1\over p}[$,
                          intervall which tends to $\R$. 
                        
                         Then 
                         $(\tilde u_\Lambda$ , $\tilde v_\Lambda)$  satisfy the  equation 
                         
                         \begin{eqnarray}\label{eq1}{|\tilde u_\Lambda^\prime |^p \over p} + {|\tilde v_\Lambda^\prime |^p \over p} &-&{\tilde u_\Lambda ^p \tilde v_\Lambda^p \over p} - \alpha {\tilde u_\Lambda^{p+2} \over (p+2)m_\Lambda^{p-2}\Lambda} -\beta {\tilde v_\Lambda ^{p+2} \over (p+2)m_\Lambda^{p-2}   \Lambda (p+2)} \nonumber \\
                         &+& \lambda_{1, \Lambda} { \tilde u_\Lambda^p \over p m_\Lambda^p \Lambda }  + \lambda_{2, \Lambda} { \tilde v_\Lambda^p \over p  m_\Lambda^p \Lambda } 
       = {T_\Lambda \over p m_\Lambda^{2p} \Lambda}.
                         \end{eqnarray}
                         
                         Using the fact that $(u_\Lambda^\prime)$ is bounded  independently on $\Lambda$, by the mean value's theorem
                         $$|\tilde u_\Lambda(y)-{1\over m_\Lambda }u(x_\Lambda) |\leq {1\over m_\Lambda^2\Lambda^{1\over p} }|u_\Lambda^\prime|_\infty$$
                         and we have analogous estimates  for $\tilde v_\Lambda$, so using $u_\Lambda (x_\Lambda) = m_\Lambda= v(x_\Lambda) $ one obtains that 
                          $\tilde u_\Lambda$ goes to $1$ uniformly,  $\tilde v_\Lambda$ goes to $1$,  finally passing to the limit in the equation (\ref{Tlambda}) one gets that 
                          $$-{1\over p} = 0,$$
                           a contradiction. 
                          We have obtained that $\Lambda m_\Lambda^{2p}$ is bounded.
                                      
                    To end the proof,  suppose by contradiction that $\Lambda m_\Lambda ^{2p} \rightarrow 0$ for 
                    a subsequence and  let $\tilde u_\Lambda(y) = {1\over m_\Lambda} u_\Lambda(m_\Lambda y+ x_\Lambda)$, 
                    then $\tilde u_\Lambda$ and $\tilde v_\Lambda$ satisfy the identity
                    
                    \begin{eqnarray*}
                     {|\tilde u_\Lambda^\prime|^p \over p} + {|\tilde v_\Lambda^\prime |^p \over p} -{m_\Lambda^{2p} \Lambda \tilde u_\Lambda^p
                    \tilde v_\Lambda^p\over p} &-&{\alpha \over( p+2)} m_\Lambda ^{p+2} \tilde u_\Lambda^{p+2} -{\beta \over (p+2)} m_\Lambda^{p+2} \tilde v_\Lambda^{p+2} \\
                    &+& {\lambda_{1,\Lambda} m_\Lambda ^p \tilde u_\Lambda^p\over p} +{ \lambda_{2, \Lambda} m_\Lambda ^p  \tilde v_\Lambda ^p\over p}= {T_\Lambda\over p}.
                    \end{eqnarray*}

                      Furthermore 
                      $|\tilde u_\Lambda(y)-{1\over m_\Lambda }u_\Lambda(x_\Lambda)|\leq C|y|$
                       by the mean values theorem, which implies that $ \tilde u_\Lambda$ is  uniformly bounded.  
                        Since $|\tilde u_\Lambda^\prime|_\infty$ is bounded independently on $\Lambda$, the dominated convergence theorem implies that  (up to subsequence) $u_\Lambda^\prime$ converges strongly in $L^p$, in particular  $|u^\prime_\Lambda|^{p-2} u_\Lambda^\prime $ converges
                         in the distributional sense, and so,  also does its derivative. 
                       By passing to the limit one obtains that  $(\tilde u_\Lambda, \tilde v_\Lambda)$ tends locally uniformly to $(u_\infty, v_\infty)$ which satisfies 
                      $| u^\prime_\infty |^{p-2} u_\infty^{\prime\prime} = 0$
                        and  $| v_\infty^\prime |^{p-2} v_\infty^{\prime\prime} = 0$, hence $u_\infty^\prime = cte$, $v_\infty ^\prime = cte$. Since   $u_\infty$ and $v_\infty$ are bounded, these constant are zero, which  yields   a contradiction with the identity  
                        $${|u_\infty^\prime|^p \over p} +{|v_\infty^\prime |^p\over p} ={ T_\infty\over p},$$
                         where we have used the fact that $T_\infty  := \lim T_\lambda \neq 0$,  by the estimates on $T_\Lambda$ proved before. We have obtained that $m_\Lambda^{2p} \Lambda $ is bounded from above by some constant $>0$.

                          We finally prove that $\Lambda^{1\over 2p} \min (x_\Lambda -a, b-x_\Lambda) \rightarrow +\infty$. Suppose  for example and by contradiction that  up to a subsequence $\Lambda^{1\over 2p}  (x_\Lambda -a)  \rightarrow  C_1$.
                          Define 
                          $\tilde u_\Lambda(y) = \Lambda^{1\over 2p} u({y\over \Lambda^{1\over 2p} } + x_\Lambda)$. 
                          Then $\tilde u_\Lambda$,  $\tilde v_\Lambda$ satisfy 
                          
                          $$|\tilde u_\Lambda^\prime|^{p-2} \tilde u_\Lambda^{\prime\prime} -{\alpha \tilde u_\Lambda^{p+1}\over \Lambda^{p+2\over 2p}} -\tilde v_\Lambda^p \tilde u_\Lambda^{p-1} +\lambda_{1, \Lambda}  {\tilde u_\Lambda^{p-1}\over \Lambda^{1\over 2}}=0,$$
                            $$|\tilde v_\Lambda^\prime|^{p-2} \tilde v_\Lambda^{\prime\prime} -{\beta \tilde v_\Lambda^{p+1}\over \Lambda^{p+2\over 2p}} -\tilde u_\Lambda^p \tilde v_\Lambda^{p-1} +\lambda_{2, \Lambda}  {\tilde v_\Lambda^{p-1}\over \Lambda^{1\over 2}}=0.$$

   We also get from the energy estimate 
   
   $${|\tilde u_\Lambda^\prime |^p\over p} + {|\tilde v_\Lambda^\prime |^p\over p}-{\tilde u_\Lambda^p \tilde v_\Lambda^p\over p} -{\alpha  \tilde u_\Lambda^{p+2}\over (p+2) \Lambda^{p+2\over 2p}} -{\beta  \tilde v_\Lambda^{p+2}\over (p+2) \Lambda^{p+2\over 2p}}+ \lambda_{1, \Lambda} {\tilde u^p \over \Lambda ^{1\over 2}}+ \lambda_{2, \Lambda} {\tilde v^p \over \Lambda ^{1\over 2}} = {T_\Lambda\over p}. $$
    Remark as before that  when $\Lambda$ goes to infinity, $(\tilde u_\Lambda, \tilde v_\Lambda)$ tends  locally uniformly to some $(U, V)$ , which satisfies 
    
    $$|U^\prime |^{p-2} U^{\prime \prime} = V^p U^{p-1}$$
     $$|V^\prime |^{p-2} V^{\prime \prime} = U^p V^{p-1}$$
and $U(-C_1) = V(-C_1) = 0$. Note that if $(b-x_\Lambda) \Lambda ^{1\over 2p}\rightarrow C_2<\infty$ one  has $U(C_2)=0$ and then $U^{\prime \prime} \geq 0$ and $U\geq 0$  implies $U\equiv 0$. We then assume that $C_2 = +\infty$.

Furthermore using Fatou's lemma one gets 

$$\int_{-C_1}^\infty V^p U^{p-1} \leq \liminf \int_{(a-x_\Lambda)\Lambda^{1\over 2p}}^{(b-x_\Lambda) \Lambda^{1\over 2p}} \tilde v^p \tilde u^{p-1} = \Lambda \int_a^b v^p u^{p-1} \leq C$$
as stated in the proof of Lemma \ref{lemLambda}. 
Since $U^\prime$ is increasing and $U \geq 0$,  if $U$ is not identically zero, there exists   $C_1^\prime \geq C_1$, such that  $U^\prime (C_1) >0$, and $U = 0$ on $[-C_1, -C_1^\prime]$  then 
$U(x) \geq U(-C_1^\prime ) + U^\prime (-C_1^\prime )(x-C_1^\prime )= U^\prime (-C_1^\prime)+ (x-C^\prime _1)$. 
 In the same manner there exists $C_1^{\prime\prime} \geq C_1$, such that 
 $$V(x) \geq V (-C_1^{\prime\prime})+ V^\prime(-C_1^{\prime\prime}) (x-C_1^{\prime\prime}). $$

 We have obtained  that near $+\infty$,  $U^{p-1} V^p \geq Cx^{2p-1}$,  which contradicts the fact that 
  $U^{p-1} V^p$ is integrable on $]-C_1,  + \infty[$.  Finally $U= V=0$,  a contradiction with the identity defining $T_\infty$, when passing to the limit.   We have obtained that 
  $\Lambda^{1\over 2} \min ( x_\Lambda-a, b-x_\Lambda)\rightarrow +\infty$.

   This ends the proof of Theorem \ref{thlambda}.

\section{Qualitative properties of the $p$-system in the one dimensional case : Proof of Theorem \ref{thUV}}

In this section  we want to prove the existence of non trivial solutions to the limit system (\ref{eq5}). Note that the previous  existence 's result  is 
obtained under the assumption $(\lambda_{1, \Lambda},\lambda_{2, \Lambda})\leq C$. 
Theorem\ref{thUV} is a consequence of several  Theorems and propositions : 
 
\begin{theo}\label{theoexi4}

There exists an entire solution  for (\ref{eq5}) such that $U(x) = V(-x)$.
\end{theo}

Proof 

We argue as in \cite{BLWZ}, up to technical arguments due to the non linearity of the  $p$-Laplacian, and due to the singularity $(p<2)$ or the degeneracy $(p>2)$. 

 Let us  consider for $R$ large the variational problem
 $$\inf_{ \{(U,V)\in H^1(]-R,R[), U(x) = V(-x), U(-R) = 0, U(R) = R\}} {1\over p} \int_{-R}^R |U^\prime|^p +   {1\over p} \int_{-R}^R |V^\prime|^p + {1\over p} \int |U|^p |V|^p. $$
  This problem admits a unique solution $(U_R, V_R)$. 
  
  We prove that $U_R$ is non negative. Indeed one has 
  $$|U_R^\prime|^{p-2} U_R^{\prime\prime} = |U_R|^{p-2} U_R |V_R|^p, $$
   and a symmetric equation for $V_R$.
   Multiplying by $U_R^-$ and integrating by parts, using the fact that $U_R(-R)$ and $U_R(R)$ are nonnegative,  one gets that $U_R^-=0$ and then $U_R\geq 0$. The same is valid for $V_R$. 
    
  By the  strong maximum principle of Vasquez  \cite{V}, $U_R >0, $ on $]-R,R[$,   $U_R^{\prime } (-R)> 0$ and $|U_R^\prime|^{p-2} U_R^\prime$ is increasing implies that $U_R^\prime>0$ everywhere, finally $U_R^{\prime\prime} \geq 0$.  
  Analogously $V_R^\prime <0$,   then $U_R-V_R $  vanishes  only on zero. 
   $U_R-V_R>0$ on $]0, R]$ implies that 
  $|U_R^\prime|^{p-2} U_R^{\prime\prime} -  |V_R^\prime|^{p-2} V_R^{\prime\prime} = U_R^{p-1} V_R^{p-1} (V_R-U_R) \leq 0$ on $[0, R]$.

    {\bf First case $p\geq 2$.}
    
            Using 
     $V_R^\prime (x) = -U_R^\prime (-x)$ one has 
      for $x>0$
      $|V_R^\prime | (x)=  U^\prime_R (-x) < U_R^\prime (x)$ and then 
       $|U_R^\prime |^{p-2} U_R^{\prime\prime} \leq |U_R^\prime|^{p-2} V_R^{\prime\prime}$ which implies since $U^\prime_R >0$, that 
       $U_R^{\prime\prime}-V_R^{\prime\prime} \leq 0$. From this one gets that  $U_R \geq V_R+ x$ for $x>0$, 
       Indeed on $0$, $U_R(0) = V_R(0) \geq V_R(0)+ 0$ and on $R$, $U_R(R) = R \geq V_R(R)+ R = R$, 
       in particular 
       $U_R(x) > x$ for $x>0$. 
    
We have obtained that 
      
  $|V_R^\prime |^{p-2} V_R^{\prime\prime} \geq  (x^+)^p V_R^{p-1}$ since $V_R^{\prime\prime} \geq 0$ anywhere else,  in particular on $\R^-$.  
  Let $\bar V$ be the solution, (  given by Lemma (\ref{smallsolutions})),  of 
  $$|\bar V^\prime|^{p-2} \bar V^{\prime\prime} = x^p |\bar V|^{p-2} \bar V,  $$ 
   on $\R^+$, which  is positive and satisfies 
  $\bar V^\prime (0)=-2$.  Let us extend $\bar V$ on $\R^-$  by the  linear function $-2x+ \bar V(0)$.  
   Since $\bar V$ hence defined is ${\cal C}^2$ and is a solution  of   $$|\bar V^\prime|^{p-2} \bar  V^{\prime\prime} 
  = (x^+)^p |\bar  V|^{p-2} \tilde V  $$ on  both $\R^+$ and $\R^-$,  one gets that it is a solution on $]-R, R[$. 
   For $R$ large enough, $V_R(-R) =R\leq  -2 (-R) + \bar V(0) = \bar V(-R)$, while 
   $V_R(R) =0 \leq \tilde V(R)$ since $\tilde V$ is positive. Using the comparison principle,  since $\bar V$ and $V_R$ are respectively solution and sub-solution of the same equation, one gets that 
   $V_R\leq \bar V$ on $[-R, R]$.
   Using  Harnack's  inequality, \cite{DS} one  gets that $(U_R, V_R)$ tends to a non trivial solution (since $U_R \geq x^+$), $(U, V)$ which satisfies 
   $V(x) = U(-x)$.

    {\bf Second case $p<2$: }
    
     We begin to prove that 
     $U_R^\prime (x) \geq {1\over 1+2^{1\over p-1}} $
      for $x>0$. 
Note that  for $x>0$
$$ (U_R^\prime )^{p-1}\leq |V_R^\prime |^{p-2} V_R^{\prime}(x) + 2 (U_R^\prime (0))^{p-1}$$
 which implies that 
 $U_R^\prime (x) \leq 2^{1\over p-1} U_R^\prime (0)$.
 Hence 
 $ U_R(R)= R \leq U_R(0) + 2^{1\over p-1}  RU_R^\prime (0)$
 On the other hand  when $x \in \R^-$
 $$U_R^\prime (x) \leq U_R^\prime (0)$$
  which implies y the mean value's theorem that 
  $U_R(0) \leq U_R(-R) + U_R^\prime (0) R$.
  We have obtained that 
  $$U_R^\prime (0) \geq {1\over 1+ 2^{1\over p-1}} $$
   as soon as $R$ is large enough.   
   We derive from this that on $\R^+$
   $ |V_R^\prime |^{p-2} V_R^{\prime \prime} \geq   \left({1\over 1+2^{1\over p-1}}\right)^{p}  (x^+)^p V_R^{p-1}.$
   
    We now consider the solution $\bar V$  of 
  $|\bar V^\prime |^{p-2} \bar V^{\prime\prime} =  \left({1\over 1+2^{1\over p-1}}\right)^{p}x^p \bar V^{p-1}  $ on $\R^+$, $\bar V>0$ given by  Proposition \ref{propbetagamma}  which satisfies $\bar V^\prime (0) = -2$, extended by $-2 x + \bar V(0)$ on $\R^-$. 
    
     One obtains as in the case $p\geq 2$ that 
     $V_R \leq \bar V$ and by Harnack's inequality, \cite{DS},  one gets that $(U_R, V_R)$ tends locally uniformly to $(U, V)$ which is not identically zero.

          \begin{lemme}
    
      Suppose that  $(U, V)$ satifies (\ref{eq5}). Then either $U^\prime >0$ and $V^\prime <0$ or $U^\prime <0$ and $V^\prime >0$. Furthermore there exists some constant $C$ such that 
      $|U^\prime | + |V^\prime | \leq C$.
      \end{lemme}
         
       Proof 
       
        Clearly the identity 
        \begin{equation}\label{eqetoile} |U^\prime |^p+
|V^\prime |^p -U^p V^p = T_\infty
\end{equation}
holds  for some  finite constant $T_\infty$. Since $U^{\prime\prime} \geq 0$, either $U^\prime>0$ or $U^\prime <0$ or there exists $x_1$ such that $U^\prime (x) >0$ for $x> x_1$ and conversely for $x < x_1$, and the same for $V^\prime$. 

   Suppose that $U^\prime $ and $V^\prime$ have the same sign somewhere, then if this sign is positif,  by the increasing behavior of $U^\prime$ and $V^\prime $,  it is true  also for $x$ large .

    In particular $U(+\infty)= V(+\infty) = +\infty$.

 {\bf The case $p\geq 2$}
 
 Let $\varphi = U+ V$.  For $x$ large enough,  $|\varphi^\prime |^{p-2} =( |U^\prime|+ |V^\prime |)^{p-2} \geq c (|U^\prime|^{p-2}+ |V^\prime|^{p-2})$ and then 
 $|\varphi^\prime |^{p-2} \varphi^{\prime\prime} \geq c  (|U^\prime|^{p-2}+ |V^\prime|^{p-2})(U^{\prime\prime} + V^{\prime \prime})
  \geq c(|U^\prime|^{p-2}U^{\prime\prime}+ |V^\prime|^{p-2}V^{\prime\prime})\geq 2^{p-1}(V^p+ U^p)\geq  (U+V)^p$ as soon 
  as $x$ is large enough since $U$ and $V$ go to infinity. 
  
 Then 
$ |\varphi^\prime |^{p-2} \varphi^{\prime\prime} \geq  \varphi^{p}$, for $x$ large enough. 

Multiplying by $\varphi^\prime >0$, one gets  that 
${d\over dx} (|\varphi^\prime|^{p}) \geq c{d\over dx}  (\varphi^{p+1})$. 
 Then 
 $(\varphi^\prime ) ^p \geq  c(\varphi^{p+1} + 1)$,    and since $\varphi$ tends to infinity,  for $x$ large enough,  
  $(\varphi^\prime ) ^p(x) \geq  {p\varphi^{p+1}(x)  \over 2}$.
   This implies 
  $\varphi^\prime \geq c\varphi^{p+1\over p}$, and then  for  $x$ large enough, 
  ${d\over dx} (-\varphi^{{-1\over p}}) \geq c_p$,
 which would imply that $\varphi$ becomes negative  for $x$ large enough. 
     
       {\bf The case $p<2$}
     
     Let $\varphi = \left((U^\prime)^{p-1} + (V^\prime)^{p-1}\right)^{1\over p-1}$, note that 
     $\varphi^{p-1} = (U^\prime)^{p-1} + (V^\prime)^{p-1}$, that 
     $\varphi \geq U^\prime + V^\prime$ and 
     $\varphi^p \leq  2^{{p\over p-1}-1} (U^\prime + V^\prime ) ^p $. 
     
      One has 
      \begin{eqnarray*}
       (\varphi^{p-1})^\prime  &=& {d\over dx} ((U^\prime )^{p-1} + (V^\prime )^{p-1})\\
      &=& (p-1) (U^{p-1} V^p + V^{p-1} U^p ) \\
     &\geq & (V+U)^p, 
     \end{eqnarray*}
    for $x$ large  by the behavior at infinity of $U$ and $V$. 
      
       Multiplying by $U^\prime + V^\prime >0$ one obtains  that 
       $$ \left({d\over dx} (\varphi^{p-1})\right) (U^\prime + V^\prime) \geq {1\over p+1}  {d\over dx} ((V+U)^{p+1}).
       $$
       On the other hand,  by the positivity of ${d\over dx} (\varphi^{p-1})$, and  $U^\prime + V^\prime$ one has 
        $${d\over dx} (\varphi^{p-1}) (U^\prime + V^\prime) \leq {p-1\over p}  {d\over dx} (\varphi^p)$$
        hence  integrating and using the fact that $U+V$ goes to infinity when $x$ goes to $+\infty$, one gets that there exists some constant $c_p$  such that for $x$ large enough 
         $$ (U^\prime + V^\prime)^p \geq c_p \varphi^p \geq c_p (U+V)^{p+1}$$
          We end as in the case  $p\geq 2$ and get an absurdity.

      If  the sign of $U^\prime$ and $V^\prime$  are both  negative somewhere  they are both  negative for $x<-x_1$. By considering the invariance of the equation by changing $x$ in $-x$, and reasoning as above one gets a contradiction.

     We have obtained that  up to exchanging $U$ and $V$, $U^\prime >0$ and $V^\prime <0$. 
     
     Suppose that $U^\prime\rightarrow +\infty$ somewhere,  then it  occurs at 
     $+\infty$ since $U^\prime$ is increasing,  in particular $U$ goes to $+\infty$ at $+\infty$,   and  using  (\ref{eqetoile}) so does   $UV$. 
  Then    $(V^{\prime}|V^\prime|^{p-2} )^\prime  = (p-1) (UV)^{p-1} U \rightarrow +\infty$,  which implies that $V^\prime$ goes to  $+\infty$ at $+\infty$, 
  a contradiction with $V^\prime <0$.  We have obtained that $U^\prime$ is bounded. 
      If $V^\prime \rightarrow-\infty$ somewhere, it occurs  at $-\infty$, then $V$ goes to $+\infty$ at $-\infty$ , by   ( \ref{eqetoile}) $UV$
       goes to $+\infty$   at $-\infty$
       and $ |U^\prime |^{p-2} U^{\prime\prime} = (UV)^{p-1} V \rightarrow +\infty$ at $-\infty$,  hence $U^\prime$ becomes $<0$ for $x$ large negative, a contradiction. 
       
       We have obtained that  $|U^\prime|+ |V^\prime|\leq C$.

     \begin{prop}\label{Tinfini}
     
      let $(U, V)$ be a solution  such that $U^\prime >0$ and $V^\prime <0$. Then 
      $V^p U^{p-1}\rightarrow 0, \ U^p V^{p-1} \rightarrow 0$ at $\pm \infty$. Furthermore  the following assertions   hold : 
      
      \begin{eqnarray}\label{eq1}
      U(-\infty)=0, U^\prime (-\infty) = 0, U^\prime (+\infty) = (T_\infty)^{1\over p} , \\
        V(+\infty)=0, V^\prime (+\infty) = 0, V^\prime (-\infty) =- (T_\infty)^{1\over p} 
        \end{eqnarray}

\end{prop}
     
 Proof 
 
  Since ${d\over dx} (|U^{\prime}|^{p-2} U^{\prime}) \geq 0$, and  ${d\over dx} (|V^{\prime}|^{p-2} V^{\prime}) \geq 0$,  $U^\prime$ and $V^\prime $ have a limit at infinity. 
Furthermore $V^\prime \leq 0$ and is increasing so it converges at $+\infty$. Its limit must be zero since if not for $x$ large enough,  $V^\prime\leq -m <0$ and  $V$ would become negative for $x$ large. 
  
 By Lemma  1.2,   $U^\prime$ is bounded.  Furthermore  it has a  positive  finite limit at infinity,
  and by  ( \ref{eqetoile}) so does $U^p V^p$. 
 Then $U$ goes to infinity,  more precisely  $U$ behaves like an increasing linear function. 
     Then $U^{p-1} V^p = {1\over U} (U^pV^p) \rightarrow 0$.  
   Furthermore  for $x> x_o$  and  $x_o$ large, 
   $|V^\prime |^{p-2} V^{\prime\prime} -V^{p-1} \geq 0$. 
   
    Let us consider      $W =V(x_o) e^{-x+ x_o}$ which satisfies 
    $|W^\prime |^{p-2} W^{\prime\prime} -W^{p-1}\leq 0$ . 
     Using  lemma  \ref{lemW} in the appendix one gets that $V\leq W$ on $[x_o, \infty[$ and then 
     $\lim_{x\rightarrow +\infty} UV = 0$ as well as 
     $\lim_{x\rightarrow +\infty}  U^p V^{p-1} = 0$. Since $V^\prime \rightarrow 0$ at infinity,
      ( \ref{eqetoile}) implies  that 
     $|U^\prime| ^p (+\infty )  = T_\infty$. In particular $T_\infty >0$. 
     
     A symmetric result holds near $-\infty$ exchanging $U$ and $V$.

      \begin{lemme}
       Let $U$ and $V$ be  as in Proposition \ref{Tinfini}. 
        There exist  some positive constants $m, M, k, K$  which depend on $T_\infty$, such that 
     $$ m e^{-Kx^2} \leq U(x)\leq Me^{-k x^2} $$
      and  some positive constants $c$, $C$ such that 
      $c|x| U(x)\leq U^\prime(x) \leq C |x|U(x)$,  for $x$ large negative,  
       and analogous inequalities for $v$ near $+\infty$.

       Furthermore  $U$ has two asymptotic lines $y=0$ at $-\infty$ and 
        $y = (T_\infty)^{1\over p}  x+ b_1$ for some $b_1\in \R$,  at $+\infty$. Similarly 
         $V$ has asymptotic lines  $0$ at $+\infty$ and 
     $y = -(T_\infty)^{1\over p} x + b_2$  for some $b_2$ at $-\infty$.

    \end{lemme}
    \begin{rema}
     The constant $k$ and $K$ can be explicitely determined in function of $T_\infty$.  
     \end{rema}
    Proof  :

      Let $k$ be such that 
      $V(x) \geq 2k|x|$ near $-\infty$, define 
      $W(x) = C e^{-k |x|^2}$ which satisfies 
      $|W^\prime |^{p-2} W^{\prime \prime } \leq  V^pW^{p-1}$ for $x$ large enough negative, 
       where the constant has been chosen in order that 
       $W (-M)= U (-M)$, and $M$ is large enough $>0$. By Lemma \ref{lemW},  using the fact that $W^\prime$  and $U^\prime$ are bounded, one gets  that 
       $U \leq W $.

        The lower bound can be obtained by considering some $K$ such that near $-\infty$
        $V(x) \leq K|x|$. 
          Consider  
         $W(x) = C e^{-Kx^2 }$ where $C$ is chosen so that 
         $W(-M) = U (-M)$ and $M$ is large.          $|W^\prime |^{p-2} W^{\prime \prime} \geq  V^p  W^{p-1}$ as soon as $x$ is large enough negative.   
          From this one derives that 
          $W\leq U$. 
          
                          We prove the assertions concerning $U$ and $U^\prime$. 
                          We begin to prove  that $U^\prime \geq U$ for $x$ large negative. Indeed let us observe that for $|x|$ large negative 
                          $|U^\prime |^{p-2} U^{\prime \prime} \geq U^{p-1}$ and then multiplying by $U^\prime$ and integrating, using $U(-\infty) = U^\prime (-\infty)=0$, one obtains 
                          ${(U^\prime)^p} \geq {U^p}$.

  To prove a better estimate, observe that    by the  behavior of
     $V$ at $-\infty$     
      $|U^\prime |^{p-2} U^{\prime\prime} \geq C |x|^p U^{p-1}$. 
       
        We  multiply by $U^\prime$ and prove that 
        $\int_{-\infty}^{-x} |t|^p U^{p-1} (t)U^\prime(t) dt \geq C |x|^p U^p$
        
         Indeed 
         $\int_{-\infty}^{-x} |t|^p U^{p-1} (t) U^\prime(t) dt = [|t|^p {U^p\over p}|_{-\infty}^{-x}+\int_{-\infty}^{-x} |t|^{p-1}U^p(t) dt 
         = |x|^p {U^p\over p}(-x) + \int_{-\infty}^x |t|^{p-1} U^{p-1}(t) U^\prime(t) dt  \geq   |x|^p {U^p\over p}(-x)$, 
          from this one  yields 
          $$ \int_{-\infty}^{-x} |t|^p U^{p-1}(t) U^\prime(t) dt \geq  |x|^p {U^p\over p}(-x) , $$
           and then  using the fact that near $-\infty$,  $U^\prime $ and $|x|^p U^p$ tend to zero,  one gets that 
           $(U^\prime)^p -C |x|^p U^p \geq 0$. 
           
            In the same manner by the behaviour of $V$ at $-\infty$ there exists $C$ such that 
             $|U^\prime |^{p-2} U^{\prime\prime} \leq C |x|^p U^{p-1}$, here we use 
              $\int_{-\infty}^{-x} |t|^p U^{p-1} (t)U^\prime(t) dt = [|t|^p {U^p\over p}|_{-\infty}^{-x}+\int_{-\infty}^{-x} |t|^{p-1}U^p(t) dt \leq   |x|^p {U^p\over p}(-x) + {1\over |x|} \int_{-\infty}^{-x} |t|^p U^{p-1}(t) U^\prime(t) dt$.
               From this as soon as $|x|>2$,  
               $$\int_{-\infty}^{-x} |t|^p U^{p-1}(t) U^\prime(t) dt \leq 2  |x|^p {U^p\over p}(-x) .$$
                We have obtained the estimate on the right 
                $U^\prime \leq C U| x|$. 
    
   To deduce from this the asymptotic of $U$ and $V$,  we use the previous estimates on $U$ for $V$ near $+\infty$. 
     So we have 
     $|V^\prime |\geq C |x| e^{-kx^2}$
      and then  by   ( \ref{eqetoile}) one  derives that  
        $(U^\prime)^p-T_\infty \geq -C |x|^pe^{-pk x^2} $, hence 
        $U^\prime(x) -T_\infty ^{1\over p}\geq    -C |x|^pe^{-kp x^2}$
         which implies by integrating that 
         $U(x)-T_\infty^{1\over p}x \geq -C$ for $x$ large positive. 
         
                 Since $U$ is convex, $(U^\prime -T_\infty^{1\over p} )$ is increasing and since it tends to $0$ at infinity, 
         it is negative, hence 
         $U-T_\infty^{1\over p} x $ is decreasing. Finally $U-(T_\infty)^{1\over p} x$ is decreasing and minorated, hence defining    
         $ b_1 = \lim_{x\rightarrow +\infty}   (U-T_\infty^{1\over p}  x)$, one has  
        $(U-T_\infty^{1\over p}  x)\geq b_1$. 
        Of course  a  symmetric result holds for $V$.

        \begin{prop}
         Let $(U, V)$ be a solution of (\ref{eq5}), with  $U(0) = V(0)=1$. Then 
         $V(y) = U(-y)$. 
         
         \end{prop}
          We assume that $U^\prime >0$, hence we are in the hypothesis of the previous propositions. 
          We can assume that $b_1 \geq b_2$, since if not one can replace $(U(x), V(x))$ by $(V(-x), U(-x))$ which exchanges $b_1$ and $b_2$.   We use the  sliding method of Beresticky and Nirenberg, \cite{BN}.

            Let  $I_\lambda = \{ x, x> \lambda\}$ and 
            $U_\lambda (x)   = U(2\lambda -x)$, $V_\lambda (x) = V(2\lambda -x)$
             Let 
             $w_1 = U-V_\lambda$, $w_2 = U_\lambda - V$. We prove in what follows that for $\lambda$ large enough and $x\in I_\lambda$, 
             $w_1(x) > 0$ as well as $w_2 > 0$.  
             From the asymptotic behaviour of $U$ and $V$,  and since $U$ is convex, 
              $U (x) \geq T_\infty^{1\over p} x + b_1$,  and by  the asymptotic behavior  of $V$ there exists $K$ such that
              $V(x)\leq -T_\infty^{1\over p}  x^-+ K$,   this implies that 
              $w_1(x) \geq T_\infty^{1\over p} (x+ (2\lambda-x)^-)+ b_1-K$ for $x\in I_\lambda$.  So by taking 
              $\lambda$ such that 
               $\lambda  T_\infty^{1\over p}  > K-b_1$ one gets that if $x\in ]\lambda, 2\lambda[$, 
                $w_1(x) \geq T_\infty^{1\over p} (\lambda  +  0)+ b_1-K>0$  and for $x> 2\lambda$, 
                $w_1(x) \geq  T_\infty^{1\over p} (x+ x-2\lambda)+ b_1-K\geq 2 \lambda T_\infty^{1\over p} + b_1-K>0 $. 
                
                 We now derive from this that $w_2$ is also $>0$ in the same $I_\lambda$ for large values of $\lambda$.  Indeed, 
                we  have 
                $$
                  |U_\lambda^\prime |^{p-2} U_\lambda^{\prime\prime} - |V^\prime|^{p-2} V^{\prime\prime} 
                  = U^p(U_\lambda^{p-1} -V^{p-1}) + U_\lambda^{p-1} (V_\lambda^p-U^p).$$
                   Multiplying this  by   $w_2^-$,  integrating   between $\lambda$ and $x$ and using 
                   $(U_\lambda-V)(\lambda) = (U-V)(\lambda) = w_1(\lambda) >0$ one gets 
                    
                    \begin{eqnarray*}
                   - \int_\lambda^x  ( |U_\lambda^\prime |^{p-2} U_\lambda^{\prime} - |V^\prime|^{p-2} V^{\prime} )(w_2^-)^\prime 
                    &+&[ {1\over p-1} (   |U_\lambda^\prime |^{p-2} U_\lambda^{\prime} - |V^\prime|^{p-2} V^{\prime} )(w_2^-))]_\lambda ^x \\
                    &=& \int_\lambda^x  U^p(U_\lambda^{p-1} -V^{p-1}) w_2^-\\
                    &+& \int_\lambda ^x  U_\lambda^{p-1} (V_\lambda^p-U^p) w_2^-\\
                    &\leq & 0
                    \end{eqnarray*}
                     Using
                     $- \int_\lambda^x  ( |U_\lambda^\prime |^{p-2} U_\lambda^{\prime\prime}- |V^\prime|^{p-2} V^{\prime\prime} )(w_2^-)^\prime\geq 0$,  $w_2^-(\lambda)=0
$ and $w_2(\infty) =0$,  as well as the  fact that $U^\prime$ and $V^\prime$ are bounded, letting $x$ go to infinity, one gets that $w_2^-=0$ and then $w_2>0$ for $x\in I_\lambda$ and $\lambda$ large enough. 

We now define 
$\lambda^\star = \inf\{\lambda >0, w_1^\mu (x) >0 \ {\rm in} \ I_\mu,\ {\rm for\ all}  \ \mu > \lambda\}$. 
By the previous observations $w_2^\mu >0$ in $I_\mu$  for all $\mu > \lambda^\star$. 
 Since $U(0) = V(0)=1$,  one has if $\lambda < 0$, by the increasing behaviour of $U-V$, 
 $(U-V)(\lambda) <0$, which implies that $\lambda^\star \geq 0$. 
 We want to prove that $\lambda^\star = 0$. 
 
 Let us observe that by  continuity  $w_1^{\lambda^\star} \geq 0$ and $w_2^{\lambda^\star} \geq 0$  on $I_{\lambda^\star}$. By the strong maximum principle $w_1^{\lambda^\star} >0$  and   $w_2^{\lambda^\star} >0$ in $I_{\lambda^\star}$. 
 By the asymptotic behaviour,   
there exists $B_1$ such that for $x< B_1<0$, 
$V(x)+T_\infty^{1\over p} x -b_2 < {b_1-b_2\over 4}$
 Take 
 $ A = \sup ( 2\lambda^\star -B_1, \lambda^\star)$ then 
  for $x > A$, and for $0<\lambda< \lambda^\star$, 
  
  $V(2\lambda-x)+T_\infty^{1\over p} (2\lambda -x) -b_2 < {b_1-b_2\over 2}$, 
   hence  for $x>A$ and $\lambda \in ]0, \lambda^\star[$, 
   $w_1^\lambda (x) -T_\infty^{1\over p} 2\lambda \geq {b_1-b_2\over 2}$. 
   
    We now observe that 
    $\inf_{[\lambda^\star, A]} w_1^{\lambda^\star} = m>0$, indeed $w_1^{\lambda^\star}  (\lambda^\star) = U(\lambda^\star)-V(\lambda^\star) > U(0)-V(0)$ since $U^\prime >0$, $V^\prime<0$, $U(0)= V(0)$ and  $\lambda^\star >0$.  By the uniform continuity of $V$ in a compact set, there exists $\eta<\lambda^\star $ such that for 
    $|\lambda-\lambda^\star|< \eta$, for all $x\in [\lambda^\star-\eta, A]$,  one has 
    $|V(2\lambda^\star-x)-V(2\lambda -x)   |\leq {m\over 2}$ and then 
    for $x > \lambda > \lambda^\star-\eta$ and $x< A$, 
    $ U(x)-V(2\lambda-x) \geq {m\over 2}$. 
     Finally 
     $\inf_{ [\lambda, A]} w_1^\lambda>0$ for $\lambda^\star-\eta< \lambda< \lambda^\star$, and  then 
     $w_1^\lambda >0$ on a neighborhood on the left of $\lambda^\star$. This contradicts the definition of $\lambda^\star$.
      We have obtained  $\lambda^\star = 0$ and then 
      $U(x) \geq V(-x)$ for $x\geq 0$. 
      
       Since  we have seen before that $w_1^0\geq 0$ implies $w_2^0 \geq 0$,  
       $ U(-x) \geq V(x)$ for $x>0$. We have obtained that 
       $U(x) \geq V(-x)$ for $x\in \R$. 
       
       Since $U(0) = V(0)$, $U(x)-V(-x)$ reaches its minimum at zero. This implies that 
       $(w_1^0)^\prime (0) = U^\prime (0)+ V^\prime (0) = 0$. 
      By the strong comparison principle one gets that 
        $U(x) = V(-x)$. We have obtained in the same time that $b_1 = b_2$.

       Part 3) in Theorem \ref{thUV} is contained in the

        \begin{prop}
      
         Suppose that $\phi, \psi$  are  bounded solutions of          
          $$\left\{\begin{array}{lc}
           (|U^\prime|^{p-2} \phi^\prime )^\prime = (p-1) U^{p-2} V^p \phi+p U^{p-1} V^{p-1} \psi\\
            (|V^\prime|^{p-2} \psi^\prime )^\prime = (p-1) V^{p-2} U^p \phi+p U^{p-1} V^{p-1} \phi. 
            \end{array}\right.$$
             Then there exists some constant $c$ such that 
             $(\phi, \psi) = c (U^\prime, V^\prime)$.

             \end{prop}

             Proof 
             
             For personal convenience  we  use minuscule letters $(u,v)$ in place of $(U, V)$.

              We do not distinguish the case $p>2$ or $p<2$ for the moment. 
              Let $\bar\phi$, $\bar\psi$ be defined as 
              $\phi = u^\prime \bar\phi$, $\psi = v^\prime \bar\psi$.  One has 
              \begin{eqnarray*}
              (|u^\prime |^{p-2} \phi^\prime)^\prime &=& (|u^\prime |^{p-2} u^{\prime\prime})^\prime   \bar\phi +p |u^\prime|^{p-2} u^{\prime\prime} \bar \phi^\prime \\
              &+&  |u^\prime |^{p-2} u^\prime \bar \phi^{\prime\prime}\\
         &=& pu^{p-1} v^p \bar \phi^\prime + p v^{p-1} u^{p-1} v^\prime \bar \phi + (p-1) u^{p-2} v^p u^\prime \bar \phi \\
         &+&|u^\prime |^{p-2} u^{\prime} \bar \phi^{\prime \prime}.
                      \end{eqnarray*}
               On the other hand using the equation satisfied by $\phi$ one gets 
               $$p |u^\prime |^{p-2} u^{\prime\prime} \bar\phi^\prime + |u^\prime |^{p-2} u^\prime \bar\phi^{\prime\prime} = p u^{p-1} v^{p-1} v^\prime (\bar \psi-\bar \phi).$$
                    In the same manner for $v$                      $$p |v^\prime |^{p-2} v^{\prime\prime} \bar\psi^\prime + |v^\prime |^{p-2} v^\prime \bar\psi^{\prime\prime} = p u^{p-1} v^{p-1} u^\prime (\bar \phi-\bar \psi).$$
                      
                      Multiplying the first equation by $u^\prime\bar\phi$ and the second one by $v^\prime \bar\psi$ ,                     one gets  
$$p |u^\prime |^{p-2} u^\prime  u^{\prime\prime}\bar \phi  \bar\phi^\prime + |u^\prime |^{p} \bar\phi \bar\phi^{\prime\prime}+
p |v^\prime |^{p-2} v^\prime  v^{\prime\prime}\bar \psi \bar\psi^\prime + |v^\prime |^{p} \bar\psi \bar\psi^{\prime\prime}= -p(uv)^{p-1} u^\prime v^\prime (\bar\psi-\bar \phi)^2.$$

 Let us now observe that 
                    $$  p |u^\prime |^{p-2} u^\prime  u^{\prime\prime}\bar \phi  \bar\phi^\prime + |u^\prime |^{p} \bar\phi \bar\phi^{\prime\prime} = 
                    (|u^\prime |^p \bar\phi\bar \phi^\prime )^\prime-|u^\prime |^p (\bar \phi^\prime )^2.                      
                    $$

                    {\bf Claim  :  $|u^\prime |^{p} \bar \phi \bar \phi^\prime $  and $|v^\prime |^{p} \bar\psi^\prime\bar \psi$ go to  zero at $+\infty$ and $-\infty$}
                    
                      This  claim will end the proof since then we will have 
                     \begin{eqnarray*}
                     0 = [ (|u^\prime |^p \bar\phi\bar \phi^\prime )+  (|v^\prime |^p \bar\phi\bar \psi^\prime )]_{-\infty}^\infty &=& \int_{\R}|u^\prime |^p (\bar \phi^\prime )^2+ |v^\prime |^p (\bar \psi^\prime )^2\\
                     &-&p\int_{\R} (uv)^{p-1} u^\prime v^\prime (\bar\psi-\bar \phi)^2\geq 0
                     \end{eqnarray*}
                     
                      and since $u^\prime v^\prime <0$ this will imply $\bar \phi^\prime = \bar \psi^\prime = 0$ and $\bar \phi= \bar \psi$. 
                      
                      In the sequel we prove the claim for $u$ and $\phi$. The result for $v$ and $\psi$ can be done by obvious symmetric arguments.

                     {\bf Proof of the claim for $u$ and $\bar\phi$}
                     
                      \begin{eqnarray*}
                      |u^\prime |^p \bar\phi \bar \phi^\prime &=& |u^\prime |^{p-2} (u^\prime \bar \phi) (u^\prime\bar \phi^\prime) \\
                      &=& |u^\prime |^{p-2} \phi u^\prime \left(  {\phi^\prime \over u^\prime }-{\phi u^{\prime \prime} \over (u^\prime )^2}\right)\\
                      &=& |u^\prime |^{p-2}  \phi \phi^\prime - |\phi|^2 {|u^\prime |^{p-2} u^{\prime \prime} \over u^\prime } \\
                       &=&  |u^\prime |^{p-2}  \phi \phi^\prime -\phi^2 {u^{p-1} v^p \over u^\prime}. 
                       \end{eqnarray*}

                       We consider separately the cases $+\infty$ and $-\infty$.

                     {\bf The case $+\infty$.}
                     
                       Since $u$ increases  like a linear function,  $u^\prime$ is  minorated by some positive constant and $v$ 
                       goes exponentially towards zero,  the term  $\phi^2 {u^{p-1} v^p \over u^\prime}
$ on the right goes to zero. 
                       
                       On the other hand  $|u^\prime |^{p-2} \phi^\prime$ tends to zero. 
                        Indeed, its derivative is  integrable for $x$ large by the asymptotic behavior of $u$ and $v$, and the fact that $\phi$ and $\psi$ are bounded. So it has a limit. Suppose that the limit is $l\neq 0$, then $\phi^\prime\sim {l\over T_\infty^{p-2\over p}}$, this contradicts $\phi$ bounded. 
 Finally 
                       $\lim_{x\rightarrow +\infty}  |u^\prime |^p \bar\phi \bar \phi^\prime(x)  = 0$.

                      {\bf The case $-\infty$. }
                      
                      We now distinguish the case $p\geq 2$ and $p<2$
                      
                    {\bf The case  $p\geq 2$.}   
                         
                    Let us recall that                       
                       there exist some   positive constants $M$ , $k$, $c$  such that near $-\infty$,  $u \leq Me^{-k x^2}$ , and   $u^\prime 
                       \geq c|x| u$ for $x$ large enough  negative. In particular   since $v$  is linear at infinity 
                       $|{u^{p-1} v^p \over u^\prime}|\leq c |x| ^{p-1}e^{-(p-2)k x^2}  $ .

              Using  $\phi$  bounded, one gets that  $\phi^2 {u^{p-1} v^p \over u^\prime}$ goes to zero at $-\infty$. 
                       Furthermore  $|u^\prime |^{p-2} \phi^\prime$ has a limit at $-\infty$ by the equation,   if this limit was  $\neq 0$,  this  would imply that 
                       $\phi^\prime$ goes to $\pm\infty$ and would contradict $\phi$ bounded. 
                        All this implies that 
                        $|u^\prime |^{p} \bar \phi\bar \phi^\prime$ tends to zero at $-\infty$.

                        {\bf The case $p<2$}
                        
                        This case is much more involved and requires several steps.

                     {\bf  Step 1 : There exists $t_p$ which goes to $-\infty$  such that 
 $(|u^\prime|^{p-2} \phi^\prime) (t_p)\rightarrow 0$.}

   Suppose for a while that there does not exist $t_p$ which goes to $-\infty$ such that  $|u^\prime|^{p-2} \phi^\prime (t_p)\rightarrow 0$.  Then   there exists $C>0$ such that  for all $t$ large negative 
   either 
  $(|u^\prime|^{p-2} \phi^\prime)(t) \geq C$ or 
   $(|u^\prime|^{p-2} \phi^\prime)(t)\leq -C$. One assumes that we are in the first case and will give at the end the arguments in the other case. Then $\phi^\prime>0 $   near $-\infty$,  hence $\phi$ has a finite limit since $\phi$ is bounded. We begin to prove that $\phi$ tends to zero  at $-\infty$. Suppose $\phi$  does not tend to zero,  then there exists $m>0$ such that  either $\phi> m$ or $\phi \leq -m $ for $x$ large negative. 
    In the first case  for some constant $c>0$ which can vary from one line to another, 
    $(|u^\prime |^{p-2} \phi^\prime)^\prime  \geq c u^{p-2} v^p\geq c u^{p-2} |x|^p$ since the last term  $pu^{p-1} v^{p-1} \psi$  tends  to zero. 
    Integrating between $-x$ and $-x_o$  large negative, one gets 
    for $x> x_o$, 
     $|u^\prime |^{p-2} \phi^\prime (-x_o) -  |u^\prime |^{p-2} \phi^\prime(-x)  \geq c \int_{-x}^{-x_o} u^{p-3}(t) u^\prime(t)  |t|^{p-1}  dt \geq  c |x|^{p-1} u^{p-2} (-x)$.

      Indeed 
      \begin{eqnarray*}
      \int_{-x}^{-x_o} u^{p-3} (t)u^\prime(t)| t|^{p-1}  dt&=& [{u^{p-2}| t|^{p-1}\over p-2}]_{-x}^{-x_o}+ {p-1\over p-2} \int_{-x}^{-x_o} u^{p-2} (t)|t|^{p-2} dt \\
      &\geq & {u^{p-2}(-x) x|^{p-1}\over p-2} -{p-1\over (2-p) |x|^2} \int_{-x}^{-x_o} u^{p-3}u^\prime  |t|^{p-1} dt
      \end{eqnarray*}
    which implies that $|u^\prime |^{p-2} \phi^\prime(-x) \leq C_1 -C_2 |x|^{p-1} u^{p-2} $
     and then in particular $\phi^\prime \leq -C_3 x^{2-p} x^{p-1} = -C_3 x$. 
        This contradicts $\phi$ bounded. 
         In the same manner if  $\phi \leq -m<0$ one gets           
         $|u^\prime |^{p-2} \phi^\prime (-x_o) -  |u^\prime |^{p-2} \phi^\prime(-x)  \leq -c  |x|^{p-1} u^{p-2} $
         which implies 
         $\phi^\prime \geq C_3 x$, 
          this still contradicts $\phi$ bounded. 
             So we are in the hypothesis that $\phi$ tends to zero and $|u^\prime |^{p-2} \phi^\prime \geq C>0$. Then 
$\phi (x) \geq \int_{-\infty}^{-x} C (u^\prime )^{2-p}(t) dt  \geq \int_{-\infty}^{-x} |t|^{1-p}u ^{1-p} (t) u^\prime(t) dt .  $

Observe that 
 
 $\int_{-\infty}^{-x}  |t|^{1-p} u^{1-p}(t) u^\prime(t) dt  \geq c |x|^{1-p} u^{2-p}$. 
  Indeed 
   \begin{eqnarray*}
   \int_{-\infty}^{-x}  |t|^{1-p} u^{1-p} u^\prime(t) dt &=&  {1\over 2-p}\int_{-\infty}^{-x}  |t|^{1-p} {d\over dt}( u^{2-p} ) dt \\
   &=&
    {1\over 2-p} [ |t|^{1-p} u^{2-p}]_{-\infty}^{-x}+ {1-p\over 2-p} \int_{-\infty}^{-x} (-t)^{-p}  |u|^{2-p}(t) dt\\
    &\geq & {1\over 2-p} [ |t|^{1-p} u^{2-p}]_{-\infty}^{-x} -{p-1\over (2-p) |x|^2}\int_{-\infty}^{-x}  |t|^{1-p} u^{1-p} u^\prime(t) dt .
    \end{eqnarray*}
    
    This ends the proof by taking $|x|$ large enough. 
    
    We have obtained that $\phi (x) \geq C |x|^{1-p} u^{2-p}$ and replacing in the equation satisfied by $\phi$ one gets 
    $$|u^\prime |^{p-2} \phi^\prime (-x_o)-|u^\prime |^{p-2} \phi^\prime (-x) \geq C  \int _{-x}^{-x_o} |t|^{1-p} |t|^pdt = C |x|^2 .$$
     From this one derives that 
     $(|u^\prime |^{p-2} \phi^\prime)(-x) \leq -Cx^2,$
      a contradiction with the assumption. 

The case where 
    $|u^\prime |^{p-2} \phi(-x)\leq  -C<0$ for $x$ large enough can be recovered by changing $\phi$ in $-\phi$,  noting  the fact that the previous computations do note use the sign of $\psi$.

      We have obtained 
   that there exists $t_p$ which goes to $-\infty$,  such that 
 $(|u^\prime|^{p-2} \phi^\prime) (t_p)\rightarrow 0$.

 {\bf Step 2:  
  $(u^\prime)^{p-2} \phi \phi^\prime$ and ${u^{p-1} \phi^2 v^p \over u^\prime}$ both tend to zero at $-\infty$. }
 
 We multiply the equation satisfied by $\phi$,  by $\phi$ and integrate between $t_p$ and $t_{p+1}$,   where $t_p$ is some subsequence decreasing to $-\infty$,  given by step 1. One obtains 
 $$  [ |u^\prime |^{p-2} \phi \phi^\prime ]_{t_{p+1}}^{t_p} = \int^{t_p}_{t_{p+1}} (\left(p-1) u^{p-2} v^p \phi^2 + |u^\prime |^{p-2} (\phi^\prime)^2\right)+  \int^{t_p}_{t_{p+1}}  p u^{p-1} v^{p-1} \phi\psi, $$ and since 
 $u^{p-1} v^{p-1} \phi\psi $ is absolutely  integrable by the estimates on $u$ and $v$ one gets with the positivity of 
 $(p-1)u^{p-2} v^p \phi^2+ |u^\prime |^{p-2} (\phi^\prime)^2 $ and summing on $p$ that $u^{p-2} v^p \phi^2$  and $|u^\prime |^{p-2}( \phi^\prime)^2 $ are  integrable. Finally for all $s$ and $t$ going to $-\infty$
 $[ |u^\prime |^{p-2} \phi \phi^\prime ]_s^t$ tends to zero, hence $|u^\prime |^{p-2} \phi \phi^\prime $ has a limit, and since it possesses a subsequence which tends to zero, this limit is  zero. 
 
  We now prove that 
  $\phi^2 {u^{p-1} v^p \over u^\prime}$ has a finite limit. For this it is enough to prove that its derivative is integrable at -$\infty$. 
  
   \begin{eqnarray*}
  ( \phi^2 {u^{p-1} v^p \over u^\prime})^\prime &=& (\phi^2)\left( {u^{p-1} v^p \over u^\prime}\right)^\prime + 2\phi\phi^\prime {u^{p-1} v^p \over u^\prime} \\
  &+& \left((p-1) u^{p-2} v^p + {u^{p-1} p v^{p-1} v^\prime \over u^\prime} -{u^{\prime\prime} u^{p-1} v^p\over (u^\prime)^2}\right)\phi^2+ 2\phi\phi^\prime {u^{p-1} v^p \over u^\prime}\\
  &=& \left((p-1) u^{p-2} v^p + {u^{p-1} p v^{p-1} v^\prime \over u^\prime}-u^{\prime\prime} (u^\prime)^{p-2}\left( {u\over u^\prime}\right)^p{v^p\over u}\right)\phi^2+ 2\phi\phi^\prime {u^{p-1} v^p \over u^\prime}\\
  &=& \left((p-1) u^{p-2} v^p + {u^{p-1} p v^{p-1} v^\prime \over u^\prime}- u^{p-2} v^p \left({u\over u^\prime}\right)^pv^p\right)\phi^2+ 2\phi\phi^\prime {u^{p-1} v^p \over u^\prime}. 
  \end{eqnarray*}
   Each of the  first three terms above can be majorized near $-\infty$ by 
   $ C u^{p-2} v^p\phi^2$
    and then  are integrable near $-\infty$.

      For $(\phi^2)^\prime {u^{p-1} v^p \over u^\prime}$ we use Cauchy Schwarz's  inequality  as follows 
      \begin{eqnarray*}
      \left\vert\phi\phi^\prime {u^{p-1} v^p \over u^\prime}\right\vert &=& |\phi^\prime| |u^\prime |^{p-2\over 2} (u^\prime)^{-p\over 2} u^{p-1} v^p |\phi|\\
      & \leq&C |\phi^\prime| (u^\prime)^{p-2\over 2} |xu|^{-p\over 2} u^{p-1} v^p |\phi|\\
      &\leq &C|\phi^\prime| (u^\prime)^{p-2\over 2} u^{{p\over 2}-1} |x|^{p\over 2} |\phi|\\
      &\leq& C (\phi^\prime)^2 (u^\prime)^{p-2} + C v^p u^{p-2} \phi^2
      \end{eqnarray*}
      We deduce that since $\phi^2 u^{p-2} v^p$ is integrable near $-\infty$ , so is $\phi^2 u^{p-1} {v^p\over u^\prime}$ and then it tends to  zero.

   Of course we would obtain symmetric properties for $\psi$ and $v$ near $+\infty$. This ends the proof.

            \section{Appendix : Global existence  uniqueness and qualitative results  for   solutions of $|y^\prime |^{p-2} y^{\prime\prime} = x^p |y|^{p-2}y$ on $\R^+$}   
     
     We begin with  the comparison principle used in section 3. 
        
           \begin{lemme}\label{lemW}
           Suppose that $a$ is some continuous and  bounded   function such that  $a(x) >0$ for $x> x_o$. 
      Suppose that $W^{\prime}$ and $V^\prime$  are  bounded at infinity, that 
      $W(x_o) = V(x_o)$, or $W^\prime (x_o) = V^\prime (x_o)$, $\lim_{x\rightarrow +\infty} W(x) = \lim_{x\rightarrow +\infty} V(x)=0$, and 
      $|W^\prime |^{p-2}W^{\prime\prime} -a(x)|W|^{p-2} W \leq 0$ for $x> x_o$,
      
      $|V^\prime |^{p-2}V^{\prime\prime} -a(x)|V|^{p-2} V\geq 0$ . 
       Then 
       $V\leq W$ for $x> x_o$.
       \end{lemme}
         
         Proof : 
          Let us multiply the difference of the equations satisfied by $V$ and $W$,  by $(V-W)^+$ and integrate by parts,  one gets 
          $\int_{x_o}^x ( |V^\prime |^{p-2}V^{\prime} - |W^\prime |^{p-2}W^{\prime})((V-W)^+)^\prime +
           \int_{x_o}^x   a(t) (|V|^{p-2} V-|W|^{p-2} W)(V-W)^+(t) dt  \leq 0$. 
           Passing to the limit when $x$ goes to infinity and using 
           $( |V^\prime |^{p-2}V^{\prime\prime} - W^\prime |^{p-2}W^{\prime\prime})((V-W)^+)^\prime>0$ one gets in particular  that 
           $ \int_{x_o}^\infty   a(t) (|V|^{p-2} V-|W|^{p-2} W)(V-W)^+(t) dt =0$ and then $V \leq W$ on $[x_o, \infty[$. 
           
          \begin{prop}\label{propexilocal}
          For $(y_o, y_1)$ given  there  exists a  unique global  solution on $\R^+$ of 
           $$ |y^\prime |^{p-2} y^{\prime\prime} = x^p |y|^{p-2} y, \ y(0) = y_o, \ y^\prime (0) = y_1.$$
          \end{prop}
           
            Proof  of  Proposition \ref{propexilocal}. 
            
            We begin to prove local existence and uniqueness of solutions.  Suppose that $x_o\geq  0$. 
             
              Let $y_o= y(x_o)$, $y_1 = y^\prime (x_o)$. If $y^\prime (x_o) \neq 0$,  Cauchy Lipschitz theorem
 can be applied and provides local existence and uniqueness of the solution. 
 
  Suppose that $y_o = y_1=0$.  Then we use some strict maximum principle to get that $y\equiv 0$ on the right and the left of $x_o$. 
  
   Suppose indeed that $y$ is not identically zero. 
    We begin to prove that  if $y(x_o+ h) >0$ for some $h>0$,  then 
    $y\geq 0$ on $[x_o, x_o+ h]$.   We multiply the equation by $y^-$ and integrate between $x_o$ and 
$x_o+h$, we get 
$\int_{x_o}^{x_o+h} |y^\prime|^{p-2} y^\prime (- (y^-)^\prime) -p \int_{x_o}^{x_o+h} x^p |y|^{p-2} y y^-+ [   |y^\prime|^{p-2} y^\prime  (y^-)]_{x_o}^{x_o+ h} =0$
 and since $y^-=0$ on $x_o$ and $x_o+ h$ one obtains 
 $\int_{x_o}^{x_o+h} |(y^-)^\prime|^{p}   + p \int_{x_o}^{x_o+h} x^p |y^-|^{p} =0$
  and then 
  $y^-=0$ on $[x_o, x_o+ h]$. 
 So we are  in the situation where $y\geq 0$ on $[x_o, x_o+h]$. We prove that if $y$ is not identically zero on the right,
  there is a contradiction with $y^\prime (x_o)=0$. 
 Let $\gamma$ 
 be such that $\gamma > x$ on $[x_o, x_o+\delta]$, $\beta $ such that $\beta (e^{\gamma \delta}-1) < y(x_o+\delta)$,  and consider  
 $w =     \beta (e^{\gamma (x-x_o)}-1)$. Then 
 $w \leq y$ on $\{x_o\}$ and on $\{x = x_o+\delta\}$, and 
 $|w^\prime |^{p-2} w^{\prime\prime } > x^p |w|^{p-2} w$. By the classical comparison principle one gets that 
 $ u \geq w$ on $[x_o, x_o+\delta]$,  which implies that 
 $$\liminf {y(x_o+ h)- y(x_o)\over h} \geq\liminf  {w(x_o+ h)-w(x_o)\over h} = \gamma \beta  >0$$ and  contradicts $y^\prime (x_o)=0$. 
  Doing the same on the left, one gets that $y\equiv 0$. 
  
   We now suppose that 
   $y_1=0$ and $y_o\neq 0$. 
   
    We use the fixed point theorem to obtain existence and uniqueness of solution.  One can suppose without loss of 
    generality that $y_o = 1$.

     In the following we suppose first that $x_o\neq 0$, and will give at the end the changes to bring when $x_o=0$, considering only the right hand side. 
    
     Define  for $y\in  B_{x_o, \delta} (1, {1\over 2}) := \{ y\in {\cal C}  ([x_o-\delta, x_o+ \delta ]), 
     \ |y-1|_{  {\cal C}  ([x_o-\delta,  x_o+ \delta ])} \leq {1\over 2}\}$. Let us  define    the function $\phi$ as  $\phi (Z) = |Z|^{{1\over p-1}-1} Z$, and the operator   $T$  as 
     $T(y)(x) = 1+ \int_{x_o} ^x \phi (\int_{x_o}^ t (p-1) s^p |y|^{p-2} y(s) ds ) dt$.
      We prove in what follows  that one can choose $\delta$ small enough in order that 
      $T$ sends  $ B_{x_o, \delta} (1, {1\over 2}) $,  into itself and is contracting  in that ball. 
      
      Indeed we use  for $y$ in  $B_{x_o, \delta} (1, {1\over 2}) $
      $(p-1) (|x_o-\delta|) ^p |t-x_o| \left({1\over 2}\right)^{p-1} \leq| \int_{x_o}^ t (p-1)  s^p |y|^{p-2} y(s) ds | \leq (p-1) (|x_o|+\delta)^p  |t-x_o| \left({3\over 2}\right)^{p-1}$. Using the mean value's theorem, denoting 
     $ Y_i(t)= (\int_{x_o}^ t (p-1)  \ s^p |y_i|^{p-2} y_i(s) ds )  $ for $i=1,2$,   one gets that for some $C$ independent on $\delta$  
      $$ |\phi (Y_1(t)) -\phi(Y_2(t) ) | \leq |Y_1(t)-Y_2(t) | (C |t-x_o|)^{{1\over p-1}-1}.$$
       We now observe that 
       $|Y_1(t)-Y_2(t) |\leq  (p-1)(|x_o|+ \delta)^p  \left({3\over 2}\right)^{p-2} ||y_1-y_2||_\infty |t-x_o|$, from this one derives that 
       
      $$|\phi(Y_1(t))-\phi(Y_2(t))| \leq  C||y_1-y_2||_\infty |t-x_o|^{1\over p-1}, $$
       and then  for $x> x_o$
       $$\int_{x_o}^x |\phi(Y_1(t))-\phi(Y_2(t))| dt \leq  C||y_1-y_2||_\infty  |x-x_o|^{p\over p-1}$$
     In particular choosing  $\delta$ such that $C \delta ^{p\over p-1}< {1\over 2
     }$,  the map $T$ is contracting.  Under the same condition $T$ maps $B_{x_o, \delta} (1, {1\over 2}) $ into itself. 
     
     Then it possesses a unique  fixed point. Since   any solution is around $x_o$ a fixed point of $T$ we have otained the local existence and uniqueness. 
   
   We now give the changes to bring when $x_o=0$. 
   
   Define  for $y\in  B_{0, \delta} (1, {1\over 2}) := \{ y\in {\cal C}  ([0,  \delta ]), 
     \ |y-1|_{   {\cal C}  [0,  \delta ])} \leq {1\over 2}\}$, 
     $T(y) = 1+ \int_0 ^x \phi (\int_{0}^ t (p-1) s^p |y|^{p-2} y(s) ds ) dt$. 
      We prove  in what follows that one can choose $\delta$ small enough in order that 
      $T$ send  $ B_{0, \delta} (1, {1\over 2}) $,  into itself and is contracting  in this ball. 
      
      Indeed we use  for $y$ in that ball 
      ${p-1\over p+1}\left( {1\over 2}\right)^{p-1}  t^{p+1} \leq| \int_{0}^ t (p-1)  s^p |y|^{p-2} y(s) ds | \leq {p-1\over p+1} t^{p+1} 
       \left({3\over 2}\right)^{p-1}$ and by  the mean value's theorem, denoting 
     $ (\int_{0}^ t (p-1)  \ s^p |y_i|^{p-2} y_i(s) ds ) = Y_i(t)$ , one gets  
      \begin{eqnarray*}
      |\phi (\int_{0}^ t (p-1) s^p |y_1|^{p-2} y_1(s) ds ) &&-\phi((\int_{0}^ t (p-1) s^p |y_2|^{p-2} y_2(s) ds ) | \\
      &\leq& |Y_1(t)-Y_2(t) | C t^{{p+1\over p-1}-(p+1)}
      \end{eqnarray*}
     We now observe that 
       $|Y_1(t)-Y_2(t) |\leq  (p-1)(|x_o|+ \delta)^p  \left({3\over 2}\right)^{p-1} ||y_1-y_2||_\infty |t|^{p+1}$, from this one derives that 
       
      $$|\phi(Y_1(t))-\phi(Y_2(t))| \leq  C||y_1-y_2||_\infty |t|^{p+1\over p-1}$$
       and then  for $x> 0$
       $$\int_{0}^x |\phi(Y_1(t))-\phi(Y_2(t))| dt \leq  C||y_1-y_2||_\infty  |x|^{2p\over p-1}.$$
        In   particular choosing  $\delta$ such that $C \delta ^{2p\over p-1}< {1\over 2
     }$,  the map $T$ is contracting.  Under the same condition $T$ maps $B_{0, \delta} (1, {1\over 2}) $ into itself.

     We want to prove global existence. 
     For that aim, suppose that there exists $\bar x$ such that either $y^\prime (\bar x) = +\infty$ or $y(\bar x) = +\infty$. 
     If $y(\bar x) < \infty$, $|y^\prime|^{p-2} y^\prime  (\bar x) = |y^\prime|^{p-2} y^\prime  (\bar x-h) + \int_{\bar x-h}^{\bar x} p t^p |y|^{p-2} y (t) dt $
      then 
      $ \lim_{x\rightarrow \bar x} y^\prime (x) $ is finite and we get local existence after $\bar x$, so we can assume that 
      $y(\bar x) =+\infty$. We begin to observe that by continuity $y>0$ in a neighborhood on the left of $\bar x$,   and 
      then by the equation $|y^\prime|^{p-2} y^\prime$  is increasing on the left of $\bar x$,   hence has a limit for  $x\rightarrow \bar x$, $x< \bar x$. Suppose for a while that this limit is $ L \leq 0$. Then one would have for  $x< \bar x$
      $$|y^\prime|^{p-2} y^\prime (x) = L +(p-1)  \int_{\bar x} ^x
t^p y^{p-1} dt <0. $$
      Then $y(x) > y(\bar x)=+\infty$ a contradiction.  We have obtained that 
      $\lim_{ x\rightarrow \bar x} y^\prime (x) = L>0$. 
    
       We now write  using the equation and the increasing behaviour of $y$ 
        on  $[x_o, \bar x]$, 
       $(y^\prime )^{p-1} (x) \leq (y^\prime)^{p-1}(x_o)+(p-1) \int_{x_o}^x t^p y^{p-1} (t) dt
        \leq     (y^\prime)^{p-1}(x_o)+(p-1)y^{p-1}(x)(\bar x)^{p+1} $.  This implies that 
       $y^\prime (x) \leq\sup (2^{2-p\over p-1}, 1) ( y^\prime (x_o) +(p-1)  (\bar x)^{p+1\over p-1} y(x))$, 
        and by integrating between $x_o$ and any $x< \bar x$, 
        $\int_{x_o}^{ x}  {  y^\prime (t) dt  \over y^\prime (x_o) +(p-1) (\bar x)^{p+1\over p-1} y(t)} \leq  \sup (2^{2-p\over p-1}, 1) (x-x_o)$, which implies that 
        $\log ( y^\prime (x_o) + (\bar x)^{p+1\over p-1} y(x))\leq C (x-x_o)$ and this contradicts the fact that $y(\bar x) =+\infty$. 
        We have obtained that $y$ is defined on  $\R^+$.

           \bigskip 
       
       We now consider the equation 
        $$ V^{\prime\prime} |V^\prime |^{p-2} = x^p |V|^{p-2} V$$
         and suppose that 
         $V(0) <0$. Then either 
         $V\leq 0$, or 
          there exists 
          $\bar x$ such that for $x> \bar x$, $V >0$. 
          
           Indeed, if we contradict this fact, there exists $\bar x_1$ which is such that $V(\bar x_1) >0$ and it is a local maximum for $V$. 
           Then $V^\prime (x_1)=0$. Since $V^\prime$ is increasing around $x_1$  by the equation,   $V^\prime (x) <0$ for $x< x_1$, $V^\prime (x) >0$ for $x > x_1$,  which contradicts  the fact that $x_1$ is a local maximum. 
           So we are in the hypothesis that  $V(x) \geq 0$ for $x$ large and by the strict maximum principle $V>0$, hence $V^\prime$ is increasing  in particular either  it  is negative
            and  in that case necessarily  tends to zero, or $V^\prime >0$ somewhere and then it remains $>0$,  
            which implies  that $V$ goes to infinity at $+\infty$.

           We want to prove that it is possible to choose $V(0)>0$ in order that for $V^\prime (0) = -2$, the solution satisfy 
           $V>0$  on $[0, \infty[$, and  $V$ and $V^\prime$, tend to zero at infinity). 
           
            \begin{lemme}\label{smallsolutions}
           There exists a solution which satisfies  on $\R^+$
             $$|y^\prime |^{p-2} y^{\prime\prime} = t^p y^{p-1}, $$
              which is positive and satisfies 
              $y^\prime (\infty) = y(\infty) = 0$. Furthermore, for $y_1<0$ given, there exists a unique positive solution   
              as above with the initial condition $y^\prime (0)= y_1$.  
              \end{lemme}

               We use the existence of a sub- and a supersolution on $[0, \infty[$ which satisfies 
                $y^\prime (\infty) = y(\infty) = 0$. 
                 Note that any positive constant is a  supersolution.  
                 
                  Let us exhibit a sub-solution.

           Let $w_2(t) = e^{ -(x^2+2x)}$. 
            Then 
            $w^\prime_2(t) = -(2x+2) \ w_2$, 
            $w_2^{\prime\prime} (t) = ((2x+2)^2-2) w_2 \geq 4x^2 w$, 
             then 
             $|w_2^\prime |^{p-2} w_2^{\prime\prime} \geq 2^p  x^p w_2^{p-1}\geq   x^p w_2^{p-1}$.

             Now we use   Perron's method on every compact set $[0, R]$, ie we define 
               $y_R= \sup \{ y,  w_2 (t) \leq y(t) \leq 1, \ {\rm on} \ [0,R], \  y \ {\rm is\ a \ sub-solution } \ {\rm on} \ [0,R], \}$.  
               
               $y_R$ is classically a solution on $[0, R]$. 
               
              The sequence $y_R$ is locally uniformly bounded and  then by Harnack's inequality, \cite{DS},  it converges  locally 
              uniformly to a solution $y$ which satisfies  
              $w_2\leq y\leq 1$. 
               Since $ y>0$ and is bounded we know by the analysis made previously that $y$  goes to zero at infinity, as well as $y^\prime$. 
               
                Let us observe that  for $V^\prime (0) <0$ given, there exists  some $V(0)>0$ such that $V$ is a solution
                 for such initial conditions,  which satisfies 
                $V(+\infty) = V^\prime (+\infty)=0$. 
                
                 Indeed, let $y$ be  the  positive solution obtained above,
                 let $V={ V^\prime (0)\over y^\prime (0)} y$, then $V$ is a positive solution  which satisfies 
                the required condition. Let us prove the uniqueness of solutions  $V$  such that $\lim V= 0$, $V^\prime$ is bounded and $V^\prime (0) $ given. Suppose for that aim that $V_i, i=1,2$ are two such solutions. 
                
                 Then substracting the equations satisfied by $V_1$ and $V_2$,  multiplying by $(V_1-V_2)$ and integrating on $[0, R]$ with $R$ large,  one gets using $(V_1-V_2)^\prime (0)=0$ and  $( |V_1^\prime |^{p-2} V_1^\prime  - |V_2^\prime |^{p-2} V_2^\prime ( V_1-V_2) (R) \rightarrow 0$ that 
                 $\int_0^R  ( |V_1^\prime |^{p-2} V_1^\prime  - |V_2^\prime |^{p-2} V_2^\prime) ( V_1^\prime -V_2^\prime ) \rightarrow 0$
                  and then $V_1 = V_2$.  
                
                                 \begin{prop}\label{propbetagamma}
                  Let $\beta , \gamma >0$ be given.  There exists a unique solution  $W$ which satisfies 
                  $$\left\{ \begin{array}{lc}
                   |W^\prime |^{p-2} W^{\prime\prime} = \beta ^p x^p W^{p-1} \ {\rm on} \ \R^+\\
                   W ^\prime (0) =- \gamma, \lim_{x\rightarrow +\infty} W(x) = 0
                   \end{array} \right.$$
                   Fiurthermore $W$ and $W^\prime$ are bounded. 
                    \end{prop}
                    
                     Proof 
                      Let $W$ be a solution of the previous equation  with $W^\prime (0) = -1$, and consider 
                      $\bar W (x) = {\gamma\over{\beta} }  W( {\beta}  x)$, then $\bar W$ satisfies 
                      \begin{eqnarray*}
                      |\bar W|^{p-2} \bar W^{\prime\prime} &= &\left( {\gamma\over {\beta}} \right)^{p-1}    {\beta} ^{p}  |W^\prime |^{p-2} W^{\prime\prime} ( \beta  x ) \\
                      &= & \beta ^p x^p \bar W^{p-1} (x)
                     \end{eqnarray*}
                        $W$ is convenient.

    \end{document}